\documentclass[12pt]{article}
\usepackage{epic,eepic}
\usepackage[dvips]{epsfig}
\usepackage{color}
\usepackage{amsmath, amscd, amssymb}
\usepackage{amssymb}
\def\draw #1 by #2 (#3){
  \vbox to #2{
    \hrule width #1 height 0pt depth 0pt
    \vfill
    \special{picture #3}
    }
  }
\def\scaleddraw #1 by #2 (#3 scaled #4){{
  \dimen0=#1 \dimen1=#2
  \divide\dimen0 by 1000 \multiply\dimen0 by #4
  \divide\dimen1 by 1000 \multiply\dimen1 by #4
  \draw \dimen0 by \dimen1 (#3 scaled #4)}
  }

\begin{document}
\renewcommand{\labelenumi}{\theenumi}
\newcommand{\qed}{\mbox{\raisebox{0.7ex}{\fbox{}}}}
\newtheorem{theorem}{Theorem}
\newtheorem{example}{Example}
\newtheorem{problem}[theorem]{Problem}
\newtheorem{defin}[theorem]{Definition}
\newtheorem{lemma}[theorem]{Lemma}
\newtheorem{corollary}[theorem]{Corollary}
\newtheorem{nt}{Note}
\newtheorem{proposition}[theorem]{Proposition}
\renewcommand{\thent}{}
\newenvironment{pf}{\medskip\noindent{\textbf{Proof}:  \hspace*{-.4cm}}\enspace}{\hfill \qed \medskip \newline}
\newenvironment{defn}{\begin{defin}\em}{\end{defin}}{\vspace{-0.5cm}}
\newenvironment{lem}{\begin{lemma}\em}{\end{lemma}}{\vspace{-0.5cm}}
\newenvironment{cor}{\begin{corollary}\em}{\end{corollary}}{\vspace{-0.5cm}}
\newenvironment{thm}{\begin{theorem} \em}{\end{theorem}}{\vspace{-0.5cm}}
\newenvironment{pbm}{\begin{problem} \em}{\end{problem}}{\vspace{-0.5cm}}
\newenvironment{note}{\begin{nt} \em}{\end{nt}}{\vspace{-0.5cm}}
\newenvironment{exa}{\begin{example} \em}{\end{example}}{\vspace{-0.5cm}}
\newenvironment{pro}{\begin{proposition} \em}{\end{proposition}}{\vspace{-0.5cm}}

\setlength{\unitlength}{12pt}
\newcommand{\comb}[2]{\mbox{$\left(\!\!\begin{array}{c}
            {#1} \\[-0.5ex] {#2} \end{array}\!\!\right)$}}
\renewcommand{\labelenumi}{(\theenumi)}
\renewcommand{\b}{\beta}
\newcounter{myfig}
\newcounter{mytab}
\def\mod{\hbox{\rm mod }}
\def\scaleddraw #1 by #2 (#3 scaled #4){{
  \dimen0=#1 \dimen1=#2
  \divide\dimen0 by 1000 \multiply\dimen0 by #4
  \divide\dimen1 by 1000 \multiply\dimen1 by #4
  \draw \dimen0 by \dimen1 (#3 scaled #4)}
  }
\newcommand{\Aut}{\mbox{\rm Aut}}
\newcommand{\w}{\omega}
\def\r{\rho}
\newcommand{\DbF}{D \times^{\phi} F}
\newcommand{\autF}{{\tiny\Aut{\scriptscriptstyle(\!F\!)}}}
\def\Cay{\mbox{\rm Cay}}
\def\a{\alpha}
\newcommand{\C}[1]{\mathcal #1}
\newcommand{\B}[1]{\mathbb #1}
\newcommand{\F}[1]{\mathfrak #1}
\title{Distance-regular graphs with large $a_1$ or $c_2$}

 \author{Jack H. Koolen$^{\,\rm 1,2,}$\footnote{This work was partially supported by the Priority Research Centers Program through the National Research Foundation of Korea (NRF) funded by the Ministry of Education, Science and Technology (Grant number 2009-0094069). JHK was also partially supported by the Basic Science Research Program through the National Research
Foundation of Korea(NRF) funded by the Ministry of Education, Science and Technology (Grant number 2009-0089826).
}
 \ \ and  Jongyook Park$^{\,\rm 1}$\\
{\small {\tt koolen@postech.ac.kr} ~~
{\tt jongyook@postech.ac.kr}}\\
{\footnotesize{$^{\rm 1}$Department of Mathematics,  POSTECH, Pohang 790-785, South Korea}}\\
{\footnotesize{$^{\rm 2}$Pohang Mathematics Institute,  POSTECH, Pohang 790-785, South Korea}}}

\date{\today}

\maketitle

\begin{abstract}
In this paper, we study distance-regular graphs $\Gamma$ that have a pair of distinct vertices, say $x$ and $y$, such that the number of common neighbors of $x$ and $y$ is about half the valency of $\Gamma$. We show that if the diameter is at least three, then such a graph, besides a finite number of exceptions, is a Taylor graph, bipartite with diameter three or a line graph.

\bigskip
\noindent
 {\bf Key Words: distance-regular graphs; Taylor graphs}
\\
\noindent
 {{\bf 2000 Mathematics Subject Classification: 05E30} }
\end{abstract}
\section{Introduction}

In this paper, we study distance-regular graphs $\Gamma$ that have a pair of distinct vertices, say $x$ and $y$, such that the number of common neighbors of $x$ and $y$ is about half the valency of $\Gamma$.

To be more precise, let $\Gamma$ be a distance-regular graph with valency $k$ and diameter $D$. If $x$ and $y$ are adjacent vertices (respectively vertices at distance two), then $a_1 := a_1(x,y)$ (respectively $c_2 :=c_2(x,y)$ ) denotes the number of common neighbors of $x$ and $y$. It is known that $a_1(x,y)$ (respectively $c_2(x,y)$) does not depend on the specific pair of vertices $x$ and $y$ at distance one (respectively two).

Brouwer and Koolen \cite{bk1} showed that if $c_2 >\frac{1}{2}k$, then $D \leq 3$, and $D=3$ implies that $\Gamma$ is either bipartite or a Taylor graph. In Proposition~\ref{3}, we slightly extend this result.

In Theorem~\ref{10}, we classify the distance-regular graphs with $a_1 \geq \frac{1}{2}k-1$ and diameter at least three. Besides the distance-regular line graphs (classified by Mohar and Shawe-Taylor \cite{mst}, cf. \cite[Theorem 4.2.16]{bcn}) and the Taylor graphs, one only finds the Johnson graph $J(7,3)$ and the halved 7-cube. This is in some sense a generalization of the classification of the claw-free distance-regular graphs of Blokhuis and Brouwer \cite{bb}, as the claw-freeness condition implies $a_1 \geq \frac{1}{2}k-1$. But they also classified the claw-free connected non-complete strongly regular graphs (i.e. distance-regular graphs with diameter two). The classification of strongly regular graphs with $a_1 \geq \frac{1}{2}k-1$ seems to be hopeless, as there are infinitely many strongly regular graphs that satisfy $a_1 \geq \frac{1}{2}k-1$, but are not line graphs, for example all the Paley graphs.
A similar situation holds for the Taylor graphs. Note that the distance two graph of a non-bipartite Taylor graph is also a Taylor graph, and hence at least one of them has $a_1 \geq \frac{k-1}{2}$ and the other one has $c_2 \geq \frac{k-1}{2}$. If $a_1 = \frac{k-1}{2}$, then it is locally a conference graph, and for any conference graph $\Delta$, we have a Taylor graph which is locally $\Delta$ (see, for example \cite[Theorem 1.5.3]{bcn}). Also there are infinitely many Taylor graphs with $a_1 \geq \frac{1}{2}k$, and hence infinitely many Taylor graphs with $c_2 \geq \frac{1}{2}k$ (see, for example \cite[Theorem 1.5.3]{bcn} and \cite[Lemma 10.12.1]{gr}).

In the last section of this paper, we will discuss distance-regular graphs with $k_2 < 2k$, where $k_2$ is the number of vertices at distance two from a fixed vertex. In particular, in Theorem~\ref{6}, we show that for fixed $\varepsilon >0$, there are only finitely many distance-regular graphs with diameter at least three and $k_2 \leq (2-\varepsilon)k$, besides the polygons and the Taylor graphs. In Theorem~\ref{11}, we classify the distance-regular graphs with $k_2 \leq \frac{3}{2}k$.

\section{Definitions and preliminaries}
All the graphs considered in this paper are finite, undirected and
simple (for unexplained terminology and more details, see \cite{bcn}). Suppose
that $\Gamma$ is a connected graph with vertex set $V(\Gamma)$ and edge set $E(\Gamma)$, where $E(\Gamma)$ consists of unordered pairs of two adjacent vertices. The {\em distance} $d(x,y)$ between
any two vertices $x$ and $y$ of $\Gamma$
is the length of a shortest path connecting $x$ and $y$ in $\Gamma$. We denote $v$ as the number of vertices of $\Gamma$ and define the {\em diameter} $D$
of $\Gamma$ as the maximum distance in $\Gamma$.  For a vertex $x \in V(\Gamma)$, define $\Gamma_i(x)$ to be the set of
vertices which are at distance precisely $i$ from $x~(0\le i\le D)$. In addition, define $\Gamma_{-1}(x) = \Gamma_{D+1}(x)
:= \emptyset$. We write $\Gamma(x)$ instead of $ \Gamma_1(x)$ and define the {\em local graph} $\Delta(x)$ at a vertex $x\in V(\Gamma)$ as the subgraph induced on $\Gamma(x)$. Let $\Delta$ be a graph.  If the local graph $\Delta(x)$ is isomorphic to $\Delta$ for any vertex $x\in\Gamma(x)$, then we say $\Gamma$ is locally $\Delta$.

A connected graph $\Gamma$ with diameter $D$ is called {\em{distance-regular}} if there are integers $b_i, c_i$ $(0\leq i\leq D)$ such that for any two vertices $x, y \in V(\Gamma)$ with $d(x, y)=i$, there are precisely $c_i$ neighbors of $y$ in $\Gamma_{i-1}(x)$ and $b_i$ neighbors of $y$ in $\Gamma_{i+1}(x)$, where we define $b_D=c_0=0$. In particular, any distance-regular graph  is regular with valency $k := b_0$. Note that a  (non-complete) connected {\em strongly regular graph} is just a distance-regular graph with diameter two. We define $a_i := k-b_i-c_i$ for notational convenience.  Note that $a_i=\mid\Gamma(y)\cap\Gamma_i(x)\mid$ holds for any two vertices $x, y$ with $d(x, y)=i$ $(0\leq i\leq D).$

 For a distance-regular graph $\Gamma$ and a vertex $x\in V(\Gamma)$, we denote $k_i:=|\Gamma_i(x)|$ and $p^i_{jh}:=|\{w| w\in\Gamma_j(x)\cap\Gamma_h(y)\}|$ for any  $y\in\Gamma_i(x)$. It is easy to see that $k_i = \frac{b_0 b_1 \cdots b_{i-1}}{c_1 c_2 \cdots c_i}$ and hence it does not depend on $x$.  The numbers $a_i$, $b_{i-1}$ and $c_i$ $(1\leq i\leq D)$ are called the {\em{intersection~numbers}}, and the array $\{b_0,b_1,\cdots,b_{D-1};c_1,c_2,\cdots,c_D\}$ is called the {\em{intersection~array}} of $\Gamma$. A distance-regular graph with intersection array $\{k,\mu,1;1,\mu,k\}$ is called a {\em Taylor graph}.

 Suppose that $\Gamma$ is a distance-regular graph with valency $k \geq 2$ and diameter $D \geq 2$, and let $A_i$ be the matrix of $\Gamma$ such that the rows and the columns of $A_i$ are indexed by the vertices of $\Gamma$ and the ($x, y$)-entry is $1$ whenever $x$ and $y$ are at distance $i$ and $0$ otherwise. We will denote the adjacency matrix of $\Gamma$ as $A$ instead of $A_1$. The eigenvalues of the graph $\Gamma$ are the eigenvalues of $A$.

 The Bose-Mesner algebra $M$ for a distance-regular graph $\Gamma$ is the matrix algebra generated by the adjacency matrix $A$ of $\Gamma$. A basis of $M$ is $\{A_i~ |~ i=0,\cdots,D\}$, where $A_0=I$. The algebra $M$ has also a basis consisting of primitive idempotents $\{E_0=\frac{1}{v}J, E_1,\cdots,E_D\}$, and $E_i$ is the orthogonal projection onto the eigenspace of $\theta_i$. Note that $M$ is closed under the componentwise multiplication $\circ$. Now, let numbers $q^k_{ij}$ ($0\leq i,j,k\leq D$) be defined by  $E_i\circ E_j=\frac{1}{v}\displaystyle\sum^D_{k=0} q^k_{ij}E_k$. The numbers $q^k_{ij}$ ($0\leq i,j,k\leq D$) are called the {\em Krein parameters} of $\Gamma$ and are always non-negative by Delsarte (cf.~\cite[Theorem 2.3.2]{bcn}).\\

Some standard properties of the intersection numbers are collected in the following lemma.

\begin{lem}(\cite[Proposition 4.1.6]{bcn})\label{pre}{\ \\}
Let $\Gamma$ be a distance-regular graph with valency $k$ and diameter $D$. Then the
following holds:\\
$(i)$ $k=b_0> b_1\geq \cdots \geq b_{D-1}~;$\\
$(ii)$ $1=c_1\leq c_2\leq \cdots \leq c_{D}~;$\\
$(iii)$ $b_i\ge c_j$ \mbox{ if }$i+j\le D~.$
\end{lem}

Suppose that $\Gamma$ is a distance-regular graph with valency $k\ge 2$ and diameter $D\ge 1$. Then $\Gamma$ has exactly $D+1$ distinct eigenvalues, namely $k=\theta_0>\theta_1>\cdots>\theta_D$ (\cite[p.128]{bcn}), and  the multiplicity of
$\theta_i$ ($0\le i\le D$) is denote by  $m_i$. For an eigenvalue $\theta$ of $\Gamma$, the sequence $(u_i)_{i=0,1,...,D}$ = $(u_i(\theta))_{i=0,1,...,D}$
satisfying $u_0$ = $u_0(\theta)$ = $1$, $u_1$ = $u_1(\theta)$ = $\theta/k$, and
\begin{center}
$c_i u_{i-1} + a_i u_i + b_i u_{i+1} = \theta u_i$ $(i=2,3,\ldots,D-1)$
\end{center}
is called the {\em standard sequence} corresponding to the eigenvalue $\theta$ (\cite[p.128]{bcn}). A sign change of $(u_i)_{i=0,1,...,D}$ is a pair $(i, j)$ with $0 \leq i < j \leq D$ such that $u_iu_j < 0$ and $u_t = 0$ for $i < t < j$.\\

In this paper we say that an intersection array is {\em feasible} if it satisfies the following four conditions:\\

\vspace{-0.3cm}\hspace{-0.6cm}$(i)$ all its intersection numbers are integral;\\

\hspace{-0.6cm}$(ii)$ all the multiplicities are positive integers;\\

\hspace{-0.6cm}$(iii)$ for any $0\leq i\leq D$, $k_ia_i$ is even;\\

\hspace{-0.6cm}$(iv)$ all Krein parameters are non-negative.\\

Recall that a {\em clique} of a graph is a set of mutually adjacent vertices and that a {\em co-clique} of a graph is a set of vertices with no edges.  A clique $\mathcal{C}$ of a distance-regular graph with valency $k$, diameter $D\geq2$ and smallest eigenvalue $\theta_D$, is called {\em Delsarte clique} if $\mathcal{C}$ contains exactly $1-\frac{k}{\theta_D}$ vertices. The {\em strong product} $G\boxtimes H$ of graphs $G$ and $H$ is a graph such that the vertex set of $G \boxtimes H$ is the Cartesian product $V(G)\times V(H)$ and any two different vertices $(u,v)$ and $(u',v')$ are adjacent in $G \boxtimes H$ if and only if ($u=u'$ or $u$ is adjacent to $u'$) and ($v=v'$ or $v$ is adjacent to $v'$). For a given positive integer $s$, the {\em $s$-clique extension} of a graph $G$ is the strong product $G\boxtimes K_s$ of $G$ and $K_s$, where $K_s$ is the complete graph (or clique) of size $s$.

 A graph $\Gamma$ is called {\em  graph of order $(s,t)$} if $\Gamma$ is locally disjoint union of $t+1$ copies of $(s+1)$-cliques. Note that, if $\Gamma$ is a distance-regular graph with $c_2=1$ and valency $k$, then $\Gamma$ is a graph of order $(s,t)$ for some $s(=a_1)$ and $t$, and hence the valency $k$ is equal to $(s+1)(t+1)$.  A {\em Terwilliger graph} is a connected non-complete graph $\Gamma$ such that, for any two vertices $u, v$ at distance two, the subgraph induced on $\Gamma(u)\cap\Gamma(v)$ in $\Gamma$ is a clique of size $\mu$ (for some fixed $\mu \geq 1$). A graph $\Gamma$ is called {\em bipartite} if it has no odd cycle. (If $\Gamma$ is a distance-regular graph with diameter $D$ and bipartite, then $a_1=a_2=\ldots=a_D=0$.)

  An {\em antipodal} graph is a connected graph $\Gamma$ with diameter $D>1$ for which being at distance 0 or $D$ is an equivalence relation. If, moreover, all equivalence classes have the same size $r$, then $\Gamma$ is also called an {\em antipodal $r$-cover}. \\

Recall the following interlacing result.

\begin{theorem}\label{1}(cf.\cite{heam})   Let $m \leq n$ be two positive integers. Let
$A$ be an $n\times n$ matrix, that is similar to a (real) symmetric matrix, and let
$B$ be a principal $m \times m$ submatrix of $A$. Then, for $i=1,\ldots , m$, $$\theta_{n-m+i}(A)\leq \theta_i(B)\leq \theta_i(A)$$
holds, where $A$ has eigenvalues $\theta_1(A) \geq \theta_2(A) \geq \ldots\geq \theta_n(A)$ and B has eigenvalues
$\theta_1(B) \geq \theta_2(B) \geq \ldots \geq \theta_m(B)$.\\
\end{theorem}

For the convenience of the reader, we give a proof of the following lemma.

\begin{lem}(\cite[Lemma 3.1]{bk1})\label{8}
Let $\Gamma$ be a distance-regular graph with diameter $D$ and valency $k$. If $D\geq3$, then $b_1\geq\frac{1}{3}k+\frac{1}{3}$.
\end{lem}
\begin{pf}
 If $b_1<\frac{1}{3}k+\frac{1}{3}$, then $a_1+1>\frac{2}{3}k-\frac{1}{3}$. Let $x$ be a vertex of $\Gamma$. As $\Delta(x)$ is not a complete graph, $\Delta(x)$ has non-adjacent vertices. So, $2(a_1+1)-(c_2-1)\leq k$, and hence $c_2\geq 2(a_1+1)-k+1>\frac{1}{3}k+\frac{1}{3}>b_1$. Lemma~\ref{pre}~$(iii)$ implies that $D\leq2$ and this is a contradiction.
\end{pf}

In \cite[Proposition 5.5.1 (ii)]{bcn}, there is an error in the statement when equality occurs. Although this is fixed in \cite{b}, for the convenience of the reader, we give a proof for this.
\begin{pro}\label{0}
Let $\Gamma$ be a distance-regular graph with diameter $D\geq3$, valency $k$ and intersection number $a_1>0$. Then $a_i+a_{i+1}\geq a_1$  for $i =1, \ldots, D-1$. If $a_i + a_{i+1} = a_1$, then $i = D-1$, $a_D = 0$, $a_{D-1}=a_1$ and $b_{D-1} = 1$.
\end{pro}
\begin{pf}
Let $i \in \{1, \ldots, D-1\}$ and
let $x, y$ be a pair of vertices of $\Gamma$ at distance $i+1$. Then
$|\Gamma_i(x)\cap\Gamma_2(y)|=p^{i+1}_{i,2} =\frac{c_{i+1}(a_i+a_{i+1}-a_1)}{c_2}\geq 0$ which gives $a_i+a_{i+1}\geq a_1$. If $a_i+a_{i+1}= a_1$, then we have $i=D-1$ by \cite[Proposition 5.5.1 (i)]{bcn}, and hence $p^D_{D-1,2}=0$. This implies that each vertex of $\Gamma(y)\cap\Gamma_{D-1}(x)$ is adjacent to each vertex of
$\Gamma(y)\cap\Gamma_D(x)$.
We will show that $a_D = 0$ by way of contradiction.  Assume $a_D>0$. This implies that the complement of the local graph at $y$, say $\overline{\Delta(y)}$, is disconnected.  Let $C$ be a connected component of $\overline{\Delta(y)}$. If $C$ is a singleton, then $a_1 = k-1$, and hence $\Gamma$ is complete. So we have $|C|\geq2$. Let  $z$ and $w$ be two adjacent vertices in $C$. Then they have at most $c_2-1$ common neighbors in $\Delta(y)$, as $z$ and $w$ are not adjacent in $\Gamma$. This means that $C$ has size at least $k-c_2+1$. This means, $k=|\overline{\Delta(y)}|\geq2(k-c_2+1)$, and hence $c_2\geq\frac{1}{2}k+1$.  On the other hand,  $\overline{\Delta(y)}$ is ($k-a_1-1$)-regular, and hence $C$ has size at least $k-a_1$ and we obtain $k=|\overline{\Delta(y)}|\geq2(k-a_1)$. This means that  $a_1\geq\frac{1}{2}k$, and hence $b_1\leq\frac{1}{2}k-1 < c_2$. This contradicts $D\geq3$. Therefore $a_D=0$.
Now we show that $b_{D-1} = 1$.
Let $z$ be a vertex of $\Gamma_{D-1}(x)$. If $b_{D-1}>1$, then for any two vertices $u$ and $v$ of $\Gamma(z)\cap\Gamma_D(x)$, we have $a_2=|\Gamma(u)\cap\Gamma_2(v)|=0$ as $a_D=0$ and $p^D_{D-1,2}=0$. But, by \cite[Proposition 5.5.6]{bcn}, we have $a_2\geq{\rm min}\{b_2,c_2\}\geq1$, as $a_1>0$ and $D\geq3$. This is a contradiction. So, $b_{D-1}=1$.
\end{pf}

\section{Distance-regular graphs with large $\mathbf{a_1}$}

In this section we classify the distance-regular graphs with $a_1\geq\frac{1}{2}k-1$ and diameter $D\geq3$. First we show that if $c_2$ or $b_2$ is large, then $D=3$ and imprimitive, where imprimitive means that the graphs are either bipartite or antipodal. This generalizes \cite[Lemma 3.14]{bk1}.

\begin{pro}\label{3}
Let $\Gamma$ be a distance-regular graph with diameter $D \geq3$ and valency $k$. If $c_2>\frac{1}{2}k$ or $b_2>\frac{1}{2}k_3$, then $D=3$ and  $\Gamma$ is  either  bipartite or a Taylor graph.
\end{pro}
\begin{pf}
Assume $c_2>\frac{1}{2}k$ or $b_2>\frac{1}{2}k_3$. Then $d(y,w)\leq2$ for a fixed vertex $x$ and $y,w\in\Gamma_2(x)$. Then $p^2_{23}=0$ which in turn implies $\frac{c_3(a_2+a_3-a_1)}{c_2}=p^3_{22}=0$ and hence $a_1=a_2+a_3$. If $a_1\neq0$, then, by Proposition~\ref{0}, $D=3$ and $\Gamma$ is a Taylor graph. So, we may assume $a_1=0$. Then $a_2=a_3=0$. Furthermore, if $b_2>\frac{1}{2}k_3$ and $c_2\leq\frac{1}{2}k$, then $k_3=\frac{kb_1b_2}{c_2c_3}>\frac{k(k-1)(k_3/2)}{(k/2)c_3}=\frac{(k-1)k_3}{c_3}$. So, $c_3>k-1$ and hence $D=3$ and $\Gamma$ is bipartite. This shows that we may assume $c_2>\frac{1}{2}k$ but this again implies $D=3$ (and $\Gamma$ is bipartite) by Lemma~\ref{pre} $(iii)$ as $c_2>b_2$.
\end{pf}

Now we show that the distance-regular Terwilliger graphs with large $a_1$ and $c_2\geq2$ are known.

\begin{pro}\label{14}\
\begin{enumerate}
  \item[$(i)$] Let $\Gamma$ be a connected non-complete strongly regular graph with valency $k$. If $c_2=1$ and $k<7(a_1+1)$, then $\Gamma$ is either the pentagon or the Petersen graph.

  \item[$(ii)$] Let $\Gamma$ be a connected non-complete strongly regular Terwilliger graph with $v$ vertices and valency $k$. If $v\leq7k$ then $\Gamma$ is either the pentagon or the Petersen graph.

  \item[$(iii)$] Let $\Gamma$ be a distance-regular Terwilliger graph with $v$ vertices, valency $k$ and diameter $D$. If $k\leq(6+\frac{8}{57})(a_1+1)$ and $c_2\geq2$, then $\Gamma$ is the icosahedron, the Doro graph (see \cite[Section 12.1]{bcn}) or the Conway-Smith graph (see \cite[Section 13.2]{bcn}).
\end{enumerate}
\end{pro}
\begin{pf}
\begin{enumerate}
\item[$(i)$] Since $c_2=1$, then $a_1+1$ divides $k$, and we obtain $(t+1)(a_1+1)=k<7(a_1+1)$ with $t\in\{1,2,3,4,5\}$. Now we will show that $a_1<t$.  Let $\textbf{C}$ be the set of all $(a_1+2)$-cliques in the graph $\Gamma$. By counting the number of pairs $(x,\mathcal{C})$, where $x$ is a vertex of the graph $\Gamma$ and $\mathcal{C}$ is a clique of $\textbf{C}$ containing $x$, in two ways, we have $| V(\Gamma)|(t+1)=(a_1+2)|\textbf{C}|$.

    Suppose that $a_1\geq t$, then $|\textbf{C}|<| V(\Gamma)|$. Let $M$ be the vertex-($(a_1+2)$-clique) incidence matrix of $\Gamma$, i.e. $M$ is the $01$-matrix whose rows and columns are indexed by the vertices and $(a_1+2)$-cliques of $\Gamma$, respectively, and the $(x,\mathcal{C})$-entry of $M$ is $1$ whenever the vertex $x$ is in the clique $\mathcal{C}$ and $0$ otherwise. Then $MM^T$ is a singular matrix, as $|\textbf{C}|<| V(\Gamma)|$, and hence $-t-1$ is an eigenvalue of $\Gamma$, as $MM^T=A+(t+1)I$, where $A$ is the adjacency matrix of $\Gamma$. As $-t-1$ is an eigenvalue of $\Gamma$, by \cite[Proposition 4.4.6]{bcn}, any clique $\mathcal{C}$ of the set $\textbf{C}$ is a Delsarte clique, and hence for all $x\in V(\Gamma)$ and all $\mathcal{C}\in\textbf{C}$, there exist $y\in \mathcal{C}$ such that $d(x,y)\leq1$ by \cite[Lemma 3]{bhk}. By considering two vertices $u$ and $v$ at distance two, and the $t+1$ cliques cliques containing $u$, it follows that $u$ and $v$ have at least $t+1$ common neighbors, which is a contradiction.

    So, $0\leq a_1<t\leq5$ holds. But except for the cases $(t,a_1)=(1,0)$ and $(t,a_1)=(2,0)$, a strongly regular graph does not exist, as the multiplicity of the second largest eigenvalue is not an integer. For $(t,a_1)=(1,0)$ and $(t,a_1)=(2,0)$, we obtain the pentagon and the Petersen graph, respectively.

\item[$(ii)$] As $1+k+\frac{b_1}{c_2}k<v\leq7k$, it follows that $b_1<6c_2$. If $c_2=1$, then by $(i)$, $\Gamma$ is either the pentagon or the Petersen graph. So, we may assume that $c_2>1$. Note that $k\geq2(a_1+1)-(c_2-1)$ (as $\Gamma$ is not a complete graph), which implies that $13c_2> c_2+2b_1\geq k+1$, and hence $k<13c_2-1<50(c_2-1)$, as $c_2\geq2$. So, there are no such strongly regular Terwilliger graph with $c_2>1$ by \cite[Corollary 1.16.6]{bcn}.

\item[$(iii)$] Let $x$ be a vertex of $\Gamma$. Then the local graph $\Delta(x)$ at $x$ is an $s$-clique extension of a strongly regular
Terwilliger graph $\Sigma$ with parameters  $\bar{v}=\frac{k}{s}$, $\bar{k}=\frac{a_1-s+1}{s}$, and $\bar{c_2}=\frac{c_2-1}{s}$ by \cite[Theorem 1.16.3]{bcn}. As $\bar{c_2}\geq1$ ($c_2\geq2$), we have $c_2-1\geq s$. If $\Sigma$ is the pentagon or the Petersen graph, then $\Delta(x)=\Sigma$ and $s=1$, and hence by \cite[Theorem 1.16.5]{bcn}, we are done in this case. So we may assume that $\bar{v}>7\bar{k}$ (by $(ii)$). As $k\leq(6+\frac{8}{57})(a_1+1)$, we obtain $s>\frac{7}{57}(a_1+1)\geq\frac{1}{50}k$. Now, as $c_2-1\geq s$, it follows $k<50(c_2-1)$, and hence we are done by \cite[Corollary 1.16.6]{bcn}.
\end{enumerate}
\end{pf}

\textbf{Remark 1.} There exist generalized quadrangles of order $(q,q)$ for any prime power $q$ (see, for example \cite[p.83]{gr}). Note that the flag graph of any generalized quadrangle of order $(q,q)$ is a distance-regular graph with $k=2q$ and $c_2=1$. This shows that there are infinitely many distance-regular graphs with $a_1>\frac{1}{7}k$. See also, Theorem~\ref{10} below.\\

Before we classify the distance-regular graphs with $a_1\geq\frac{1}{2}k-1$, we first introduce some results for this classification.

\begin{lem}\label{7}
Let $\Gamma$ be a distance-regular graph with diameter $D\geq3$ and valency $k$. If $a_1\geq\frac{1}{2}k-1$ and $c_2\geq2$, then $b_2<c_2$, and hence $D=3$.
\end{lem}
\begin{pf}
If $\Gamma$ is a Terwilliger graph, then $\Gamma$ is the icosahedron by Proposition~\ref{14} $(iii)$. If $\Gamma$ is not a Terwilliger graph, then $\Gamma$ has a quadrangle. Then by \cite[Theorem 5.2.1]{bcn}, $c_2-b_2\geq c_1-b_1+a_1+2\geq 2$ holds, and hence $D=3$ by Lemma~\ref{pre} $(iii)$.
\end{pf}

\subsection{Some eigenvalues results}

In the next three lemmas, we give some results on the eigenvalues of a distance-regular graph.

\begin{lem}\label{9}
Let $\Gamma$ be a distance-regular graph with diameter three and distinct eigenvalues $k=\theta_0>\theta_1>\theta_2>\theta_3$. If $a_3=0$, then $\theta_1>0>-1\geq
\theta_2\geq-b_2\geq\theta_3$.
\end{lem}
\begin{pf}
As $a_3=0$, we know that $\theta_1$, $\theta_2$ and $\theta_3$ are the eigenvalues of $$ T:=\left[ \begin{array}{ccc}
                                                                                           -1 & b_1 & 0     \\
                                                                                            1 & k-b_1-c_2 & b_2   \\
                                                                                            0 & c_2 & -b_2  \\
                                                                                           \end{array}    \right]$$ by \cite[p.130]{bcn}.
Since the principal submatrix $\left[\begin{array}{cc}
                                     -1 & 0\\
                                      0 & -b_2\\
                                      \end{array}\right]$ of $T$ has eigenvalues $-1$ and $-b_2$, it follows that the inequality $\theta_1\geq-1\geq \theta_2\geq-b_2\geq\theta_3$ holds by
                                      Theorem~\ref{1}. As $\Gamma$ has an induced path $P$ of length three, $\theta_1\geq$ second largest eigenvalue of $P$, which is greater than zero.

\end{pf}

\begin{lem}\label{12}
Let $\Gamma$ be a distance-regular graph with diameter $D\geq3$ and distinct eigenvalues $k=\theta_0>\theta_1>\cdots>\theta_D$. If $\Gamma$ has an eigenvalue $\theta$  with multiplicity smaller than $\frac{1}{2}k$ then the following holds:\\
$(1)$ $\theta\in\{\theta_1, \theta_D\}$,\\
$(2)$ $\theta$ is integral,\\
$(3)$ $\theta+1$ divides $b_1$.
\end{lem}

\begin{pf}
 This lemma follows immediately from \cite[Theorem 4.4.4]{bcn}.
\end{pf}

\begin{lem}\label{13}
Let $\Gamma$ be a distance-regular graph with diameter $D\geq3$ and distinct eigenvalues $k=\theta_0>\theta_1>\cdots>\theta_D$. If $\Gamma$ has an eigenvalue $\theta$  with multiplicity at most $k-2$ then the following holds:\\
$(1)$ $\theta_2\geq-1-\frac{b_1}{\theta+1}$ if $\theta=\theta_1$,\\
$(2)$ $\theta_{D-1}\leq-1-\frac{b_1}{\theta+1}$ if $\theta=\theta_D$.
\end{lem}
\begin{pf}
Let $x$ be a vertex of $\Gamma$. Then the local graph $\Delta(x)$ has smallest eigenvalue at least $-1-\frac{b_1}{\theta_1+1}$ and second largest eigenvalue at most $-1-\frac{b_1}{\theta_D+1}$, by \cite[Theorem 4.4.3]{bcn}. As $\Delta(x)$ has $k$ eigenvalues at least $-1-\frac{b_1}{\theta_1+1}$ and $k-1$ eigenvalues at most $-1-\frac{b_1}{\theta_D+1}$, by Theorem~\ref{1}, the inequalities follow.
\end{pf}

Lemma~\ref{13} gives evidence for the following conjecture. \\

\textbf{Conjecture A.} Let $\Gamma$ be a distance-regular graph with diameter three and distinct eigenvalues $k=\theta_0>\theta_1>\theta_2>\theta_3$. Then $$-1-\frac{b_1}{\theta_3+1}\geq \theta_2 \geq -1-\frac{b_1}{\theta_1+1}.$$\\

\textbf{Remark 2.} Conjecture A is true when $\Gamma$ is bipartite, as then $b_1=k-1$, $\theta_1=-\theta_2=\sqrt{b_2}$ and $\theta_3=-k$. Hence $0=-1-\frac{b_1}{\theta_3+1}>\theta_2=-\sqrt{b_2}>-1-\frac{b_1}{\theta_1+1}$, where the last inequality holds, as $k>b_2$. Conjecture A is also true when $\Gamma$ is antipodal, as then $\theta_2=-1$. We also checked that all the feasible intersection arrays in the table of primitive distance-regular graphs with diameter three \cite[p. 425-431]{bcn}, satisfy this conjecture.\\

\textbf{Remark 3.} For a distance-regular graph with diameter $D\geq4$, we have $\theta_2\geq0$ by Theorem~\ref{1}.

\subsection{Classification of distance-regular graphs with $\mathbf{a_1\geq\frac{1}{2}k-1}$}

Now we are ready to classify the distance-regular graphs with $a_1\geq\frac{1}{2}k-1$.

\begin{thm}\label{10}
Let $\Gamma$ be a distance-regular graph with diameter $D\geq3$ and valency $k$. If $a_1\geq\frac{1}{2}k-1$, then one of the following holds:\\
~~~~~$(1)$ $\Gamma$ is a polygon,\\
~~~~~$(2)$ $\Gamma$ is the line graph of a Moore graph,\\
~~~~~$(3)$ $\Gamma$ is the flag graph of a regular generalized $D$-gon of order $(s,s)$ for some $s$,\\
~~~~~$(4)$ $\Gamma$ is a Taylor graph,\\
~~~~~$(5)$ $\Gamma$ is the Johnson graph $J(7,3)$, \\
~~~~~$(6)$ $\Gamma$ is the halved $7$-cube.
\end{thm}
\begin{pf}
If $\Gamma$ does not contain a $K_{1,3}$ as an induced subgraph of $\Gamma$, then one of $(1)$, $(2)$, $(3)$ or $(4)$ holds by \cite[Theorem 1.2]{bb}. So, we may assume that $\Gamma$ contains a
$K_{1,3}$. Since $\Gamma$ contains a $K_{1,3}$, we find $3(a_1+1)-3(c_2-1)\leq k$, and hence $c_2\geq\frac{1}{6}k+1$. So,  $c_2\geq2$ and it
follows that $D=3$ by Lemma~\ref{7}. We will consider two cases, namely the case $a_3=0$ and the case $a_3\neq 0$.
\begin{description}
\item[Case 1)]  First let us assume $a_3=0$. Then $c_3=k$ and it follows $k_3=\frac{kb_1b_2}{c_2c_3}=\frac{b_1b_2}{c_2}$. Let $x$ be a vertex of $\Gamma$ and $y$ be a vertex of
    $\Gamma_3(x)$. Then $\Gamma_3(x)$ contains $\frac{k(b_2-1)}{c_2}$ vertices which are at distance $2$ from $y$, as $a_3=0$. Hence $b_1b_2>k(b_2-1)$ and this implies $b_2=1$, as $b_1\leq\frac{1}{2}k$. Thus, $k_3\in\{1,2\}$, as $b_1\leq\frac{1}{2}k$ and $c_2>\frac{1}{6}k$.
 \begin{enumerate}
   \item If $k_3=1$, then $\Gamma$ is a Taylor graph.\
   \item If $k_3=2$, then $b_1=2c_2$, and hence $v=3k+3$, as $k_2=2k$. Since $b_1\leq\frac{1}{2}k$, we have $c_2\leq\frac{1}{4}k$, which implies
       $a_2\geq\frac{3}{4}k-1$, as $b_2=1$. Let $k>\theta_1>\theta_2>\theta_3$ be the distinct eigenvalues of $\Gamma$. Then
       $k+\theta_1+\theta_2+\theta_3=a_1+a_2+a_3\geq(\frac{1}{2}k-1)+(\frac{3}{4}k-1)$ implies $\theta_1+\theta_2+\theta_3\geq\frac{1}{4}k-2$. As $a_3=0$ and $b_2=1$, we have $\theta_2=-1>\theta_3$ by Lemma~\ref{9}, and hence $\theta_1>\frac{1}{4}k$.

        If the multiplicity $m_1$ of $\theta_1$ is smaller than $\frac{1}{2}k$, then by Lemma~\ref{12}, $\frac{b_1}{\theta_1+1}$ is an integer and hence $\theta_1=b_1-1$, as $b_1\leq\frac{1}{2}k$ and $\theta_1>\frac{1}{4}k$. Then
       there is no such graph by \cite[Theorem 4.4.11]{bcn}.

       If $m_1\geq\frac{1}{2}k$, then $4k^2\geq(3k+3)k=vk=k^2+m_1\theta_1^2+m_2\theta_2^2+m_3\theta_3^2>k^2+\frac{1}2{k}(\frac{1}{4}k)^2$, and hence $k<96$.

        We checked by computer the feasible intersection arrays of antipodal $3$-covers with diameter three, satisfying $k<96$, $a_1\geq\frac{1}{2}k-1$ and $c_2>\frac{1}{6}k$, and no intersection arrays were feasible .
    \end{enumerate}

\item[Case 2)] Now we assume $a_3\neq0$. We first will show that $k$ in this case $k\leq945$ holds. If $\theta_1\geq\frac{k}{4}$ holds, then similarly as in $(2)$ of \textbf{Case 1)}, we can show that $k<96$ or $\theta_1=b_1-1$. If $\theta_1=b_1-1\geq\frac{k}{4}$ holds, then $\Gamma$ is either the Johnson graph $J(7,3)$ or the halved $7$-cube by \cite[Theorem 4.4.11]{bcn}. So, we find that if $\theta_1\geq\frac{k}{4}$ holds, then we have $k<96$. So we may assume that $\theta_1<\frac{k}{4}$. As $\theta_1\geq{\rm min}\{\frac{a_1+\sqrt{a_1^2+4k}}{2},a_3\}$ (\cite[Lemma 6]{kp}) and $\frac{a_1+\sqrt{a_1^2+4k}}{2}>\frac{k}{4}$, we have $a_3<\frac{k}{4}$, and hence $c_3>\frac{3}{4}k$. As $a_3\neq 0$, $\Gamma$ is not a Taylor graph. Then by Proposition~\ref{3}, we find $\frac{\frac{k}{2}}{c_2}k\geq\frac{b_1}{c_2}k=k_2\geq2c_3>\frac{3}{2}k$, which implies $c_2<\frac{k}{3}$.  Let $\eta_1$=max$\{-1,a_3-b_2\}$ and $\eta_2$=min$\{-1,a_3-b_2\}$.
    Then $\theta_1\geq\eta_1\geq\theta_2\geq\eta_2\geq\theta_3$ by \cite[Proposition 3.2]{kpy}. Now we will show that $\theta_1\geq\frac{k}{2}-c_2$.\\
    If $\theta_2\geq0$, then $\eta_1=a_3-b_2$ and $\eta_2=-1$, and hence $a_3\geq b_2+\theta_2$ and $\theta_3\leq -1$. We find $\theta_1+\theta_2+\theta_3=a_1+a_2+a_3-k\geq
    \frac{-k}{2}-1+a_2+b_2+\theta_2$ and this implies $\theta_1\geq a_2+b_2-\frac{k}{2}=\frac{k}{2}-c_2$.\\
    If $\theta_2<0$, then $\theta_3\leq a_3-b_2$ implies $a_3\geq b_2+\theta_3$. We find $\theta_1+\theta_2+\theta_3=a_1+a_2+a_3-k\geq
    \frac{-k}{2}-1+a_2+b_2+\theta_3$ and this implies $\theta_1\geq a_2+b_2-\frac{k}{2}=\frac{k}{2}-c_2$ as $\theta_2\leq-1$ or $a_3-b_2\geq\theta_2\geq-1\geq\theta_3$. So, we have shown $\theta_1\geq\frac{k}{2}-c_2$ and this implies $\theta_1>\frac{k}{6}$ and $c_2>\frac{k}{4}$, as $c_2<\frac{k}{3}$ and $\theta_1<\frac{k}{4}$ respectively.

    If the multiplicity $m_1$ of $\theta_1$ is smaller than $\frac{1}{2}k$, then by Lemma~\ref{12}, $\frac{b_1}{\theta_1+1}$ is an integer and hence $\frac{b_1}{\theta_1+1}\in\{1,2\}$, as $b_1\leq\frac{1}{2}k$ and $\theta_1>\frac{1}{6}k$.

    If $\frac{b_1}{\theta_1+1}=1$, then there is no such graph by \cite[Theorem 4.4.11]{bcn}.

    If $\frac{b_1}{\theta_1+1}=2$, then $u_2(\theta_1)=\frac{(\theta_1-a_1)u_1(\theta_1)-1}{b_1}=\frac{\frac{3b_1}{4}-\frac{k}{2}-\frac{3}{2}}{k}<-\frac{1}{8}$, as $\theta_1-a_1=\frac{b_1}{2}-1-a_1=\frac{3}{2}b_1-k$ and $u_1(\theta_1)=\frac{\theta_1}{k}=\frac{\frac{b_1}{2}-1}{k}$. As $\theta_1=a_3$ implies
    $u_2(\theta_1)=0$, we assume $\theta_1\neq a_3$. Then $c_3u_2(\theta_1)+a_3u_3(\theta_1)=\theta_1u_3(\theta_1)$ follows that
    $u_3(\theta_1)=\frac{c_3}{\theta_1-a_3}u_2(\theta_1)$. As $u_3(\theta_1)\geq-1$, we find $\theta_1-a_3>\frac{3}{32}k$ and hence $a_3<\frac{5}{32}k$. As $\frac{b_1}{2}-1=\theta_1\geq\frac{1}{2}k-c_2$ and $k_2\geq2c_3$, we find that $\theta_1\leq\frac{9}{40}k$ implies $c_3\leq\frac{37}{44}k$, unless $k\leq160$. This shows that $k\leq160$ or $\theta_1>\frac{9}{40}k$. So, we may assume $\theta_1>\frac{9}{40}k$ and $a_3<\frac{5}{32}k$. As $a_1\geq\frac{1}{2}k-1$ and $c_2\geq2$, we have $b_2<c_2$ and this shows $k_3=k\frac{b_1b_2}{c_2c_3}<k\frac{b_1}{c_3}<k\frac{k/2}{27k/32}=\frac{16}{27}k$. Then by \cite[Theorem 4.1.4]{bcn}, we find $m_1=\frac{v}{\sum^3_{i=0}u_i(\theta_1)^2k_i}<\frac{k+k_2+16k/27}{(9/40)^2k+(-1/8)^k_2}=\frac{43k/27+k_2}{81k/1600+k_2/64}$. As $\frac{k}{3}\leq b_1\leq\frac{k}{2}$ and $\frac{k}{4}<c_2<\frac{k}{3}$, we have $k\leq k_2 < 2k$ and $\frac{43k/27+k_2}{81k/1600+k_2/64}$ has maximum value smaller than 43.88 when $k_2=2k$. So, $m_1\leq 43$. Thus, we have $k\leq \frac{(43+2)(43-1)}{2}=945$ by \cite[Theorem 5.3.2]{bcn}.

    If $m_1\geq\frac{1}{2}k$, then similarly as in $(2)$ of \textbf{Case 1)} we obtain $k<312$.

    We checked by computer the feasible intersection arrays of distance-regular graphs with diameter three satisfying $k\leq 945$, $a_1\geq\frac{1}{2}k-1$, $c_2 >\frac{1}{6}k$ and $a_3\neq0$, and only the intersection arrays of the Johnson graph $J(7,3)$ and the halved 7-cube, were feasible.
\end{description}
\end{pf}

\section{Distance-regular graphs with small $\mathbf{k_2}$}

In this section we give two results on the distance-regular graphs with small $k_2$. In the first result we look at $k_2<2k$, and in the second result we classify the distance-regular graphs with $k_2\leq\frac{3}{2}k$ and diameter at least three.

\begin{thm}\label{6}
Let $\varepsilon>0$. Then there exist a real number $\kappa=\kappa(\varepsilon)\geq3$ such that  if $\Gamma$ is a distance-regular graph with diameter $D\geq3$, valency $k\geq\kappa(\varepsilon)$ and $k_2\leq(2-\varepsilon)k$, then $D=3$ and $\Gamma$ is either bipartite or a Taylor graph.
\end{thm}
\begin{pf}
We may assume $k\geq3$. If $c_2=1$, then $b_1=1$, as $k\leq k_2=\frac{kb_1}{c_2}\leq(2-\varepsilon)k$ and this implies $a_1=k-2$. As $c_2=1$, we know that $a_1+1=k-1$ divides $k$, which gives $k=2$, a contradiction to $k\geq3$. So, we may assume $c_2\geq 2$. Suppose $\Gamma$ either has $D\geq4$ or ($D=3$ and $\Gamma$ is not bipartite or a Taylor graph), then by Proposition~\ref{3}, $c_2\leq\frac{1}{2}k$ and $b_2\leq\frac{1}{2}k_3$. As $c_2\leq\frac{1}{2}k$ (respectively $b_2\leq\frac{1}{2}k_3$), we have $(2-\varepsilon)k\geq k_2\geq2b_1$ (respectively $(2-\varepsilon)k\geq k_2\geq 2c_3$), and hence $a_1\geq \frac{\varepsilon}{2}k-1$ (respectively $a_3\geq \frac{\varepsilon}{2}k$). So, by \cite[Lemma 6]{kp}, we have $\theta_1\geq{\rm min}\{\frac{a_1+\sqrt{a_1^2+4k}}{2}, a_3\}\geq{\rm min}\{a_1+1, a_3\}\geq\frac{\varepsilon}{2}k$. This implies $u_1(\theta_1)=\frac{\theta_1}{k}\geq\frac{\varepsilon}{2}$. As $2k>(2-\varepsilon)k\geq k_2=\frac{kb_1}{c_2}$, we find $2c_2>b_1$. Then, by \cite[Lemma 5.2]{bhk2}, we obtain $D\leq 4$. For $D=3$, it is easy to check $v\leq 7k$. If $D=4$, then by \cite[Theorem 5.4.1]{bcn}, $c_3\geq\frac{3}{2}c_2$ and this implies $k_3=k_2\frac{b_2}{c_3}\leq k_2\frac{b_1}{3/2c_2}<\frac{8}{3}k$ and $k_4<\frac{32}{9}k$, as $\frac{b_1}{c_2}<2$. So, $v\leq 10k$. Then, by \cite[Theorem 4.1.4]{bcn}, the multiplicity $m_1$ of $\theta_1$ is smaller than $\frac{v}{u_1(\theta_1)^2k}<\frac{40}{\varepsilon^2}$. So, $k<\frac{({40}/{\varepsilon^2}-1)({40}/{\varepsilon^2}+2)}{2}$ (by \cite[Theorem 5.3.2]{bcn}). Thus, if we take $\kappa(\varepsilon)=\frac{({40}/{\varepsilon^2}-1)({40}/{\varepsilon^2}+2)}{2}$, then $D=3$ and $\Gamma$ is either bipartite or a Taylor graph.
\end{pf}

 \textbf{Remark 4.} The Hadamard graphs have intersection array $\{k, k-1, \frac{1}{2}k, 1; 1, \frac{1}{2}k, k-1, k\}$ and have $k_2 = 2(k-1)$. For $k = 2^t$  $(t=1,2, \ldots)$, there exists a Hadamard graph (see, for example \cite [Section 1.8]{bcn}). This shows that the above theorem is quite sharp.\\

 \textbf{Question.} For fixed positive constant $C$, are there only finitely many primitive distance-regular graphs with diameter at least three, valency $k\geq3$ and $k_2 < Ck$?\\

In the next theorem, we classify the distance-regular graphs with $k_2 \leq \frac{3}{2}k$ and diameter at least three.

\begin{thm}\label{11}
Let $\Gamma$ be a distance-regular graph with diameter $D\geq3$, $v$ vertices and valency $k$. If $k_2\leq\frac{3}{2}k$, then one of the following holds:\\
$(1)$ $k=2$ and $\Gamma$ is a polygon,\\
$(2)$ $D=3$ and $\Gamma$ is bipartite,\\
$(3)$ $D=3$ and $\Gamma$ is a Taylor graph,\\
$(4)$ $\Gamma$ is the Johnson graph $J(7,3)$,\\
$(5)$ $\Gamma$ is the $4$-cube.\\

\end{thm}

\begin{pf}
If $c_2=1$, then $b_1=1$, as $k\leq k_2=\frac{kb_1}{c_2}\leq\frac{3}{2}k$ and this implies $a_1=k-2$. As $c_2=1$, it follows that $a_1+1=k-1$ divides $k$, which gives $k=2$ and  $\Gamma$ is a polygon.  So, we may assume  $\frac{1}{2}k\geq c_2\geq2$, as $c_2>\frac{1}{2}k$ implies that $D=3$ and $\Gamma$ is either bipartite or a Taylor graph (by Proposition~\ref{3}). As the distance-regular line graphs with $c_2\geq2$ have $D=2$ (see \cite[Theorem 4.2.16]{bcn}), $c_2\geq2$ implies that  when $a_1\geq\frac{1}{2}k-1$ the graph $\Gamma$ is either a Taylor graph, the Johnson graph $J(7,3)$ or the halved 7-cube by Theorem~\ref{10}, but the halved 7-cube has $k=12$ and $k_2=35$.

 Hence, we may assume $\frac{1}{2}k\geq c_2\geq2$ and $a_1<\frac{1}{2}k-1$. This implies $b_1>\frac{1}{2}k$, and hence  $c_2\geq\frac{2}{3}b_1>\frac{1}{3}k$, as $k_2\leq \frac{3}{2}k$.
If $a_1=0$, then $\frac{1}{2}k\geq c_2\geq \frac{2}{3}(k-1)$, which implies $k\leq4$ and it is easy to check that the theorem holds in this case by \cite{bbs, bk}. So, we may assume  $a_1>0$ and this implies $k\geq5$, as $a_1<\frac{1}{2}k-1$. As $a_1>0$, we find $a_2\geq{\rm min}\{b_2,c_2\}$, by \cite[Proposition 5.5.6]{bcn}, which in turn implies $b_2<\frac{1}{3}k<c_2$, and hence $D=3$.

 So, from now on, we assume that the diameter $D$ is three, $a_1$ is positive and the valency $k$ is at least five. Here note that $v\leq \frac{7}{2}k$, as $k_3\leq k_2\frac{b_2}{c_2}<\frac{3}{2}k\frac{k/6}{k/3}=\frac{3}{4}k$ when $b_2\leq\frac{1}{6}k$, and  $b_2>\frac{1}{6}k$ implies $c_3>\frac{1}{2}k$ by \cite[Theorem 5.4.1]{bcn}, and hence $k_3\leq\frac{3}{2}k\frac{k/3}{k/2}=k$. As $k_2\geq2b_1$ and $k_2\geq 2c_3$ (Proposition~\ref{3}), we find $a_1\geq\frac{1}{4}k-1$ and $a_3\geq\frac{1}{4}k$ respectively. By \cite[Lemma 6]{kp}, we find $\theta_1\geq{\rm min\{\frac{a_1+\sqrt{a_1^2+ 4k}}{2}, a_3\}}\geq{\rm min}\{a_1+1,a_3\}\geq\frac{1}{4}k$.\vspace{0.2cm}

 If  $m_1\geq\frac{1}{2}k$, then $\frac{7}{2}k^2\geq vk=k^2+m_1{\theta_1}^2+m_2{\theta_2}^2+m_3{\theta_3}^2\geq k^2+\frac{1}{2}k(\frac{1}{4}k)^2$ implies $k\leq80$. We checked by computer the feasible intersection arrays of distance-regular graphs with diameter three satisfying $k\leq80$ and $k_2\leq\frac{3}{2}k$, and no intersection arrays were feasible.

  If $m_1<\frac{1}{2}k$, then by Lemma~\ref{12},  $\frac{b_1}{\theta_1+1}\in\{1,2\}$, as $\theta_1\geq\frac{1}{4}k$ and $b_1\leq\frac{3}{4}k$.  If $\frac{b_1}{\theta_1+1}=1$, then there is no such distance-regular graph by \cite[Theorem 4.4.11]{bcn}.

   For $\frac{b_1}{\theta_1+1}=2$, we first show that one of $m_1<48$ and $k\leq510$ holds. If $m_1\geq48$  and $k>510$, then by \cite[Theorem 4.1.4]{bcn}, $48\leq m_1=\frac{v}{\sum^3_{i=0}u_i(\theta_1)^2k_i}<\frac{7k/2}{(\theta_1/k)^2k}$, as $u_1(\theta_1)=\frac{\theta_1}{k}$ and $v\leq\frac{7}{2}k$. This implies $(\frac{\theta_1}{k})^2<\frac{7}{96}$ and hence $\theta_1<(0.271)k$. As $\frac{b_1}{2}-1=\theta_1<(0.271)k$, we have $b_1<(0.542)k+2$, which implies $a_1>(0.458)k-3$, and hence $a_1+1>(0.458)k-3>(0.271)k>\theta_1$, as $k>510$. This in turn implies $\theta_1\geq a_3$, as $\theta_1\geq{\rm min}\{a_1+1,a_3\}$.  Note that if $k_3<\frac{1}{2}k$, then $v\leq 3k$, as $k_2\leq\frac{3}{2}k$. As $\theta_1\geq\frac{1}{4}k$, we find $m_1=\frac{v}{\sum^3_{i=0}u_i(\theta_1)^2k_i}<\frac{3k}{k/16}=48$, and this contradicts $m_1\geq48$. So, $k_3\geq\frac{1}{2}k$ and this implies $b_2\geq\frac{1}{3}c_3$, as $k_3=k_2\frac{b_2}{c_3}$ and $k_2\leq\frac{3}{2}k$. Then, as $a_3\leq\theta_1<(0.271)k$, we find $c_3>(0.729)k$, and hence $b_2>(0.243)k$. Since $a_2=k-b_2-c_2$ and $c_2>\frac{1}{3}k>(0.333)k$, we find $a_2<(0.424)k$ and this implies $(0.924)k-1>a_1+a_2\geq k+\theta_2+\theta_3\geq k-3+\theta_3$, as $\theta_1\geq a_3$, and by Lemma~\ref{13}, we have $\theta_2\geq -3$. So, we obtain $-(0.076)k+2\geq \theta_3$.  As $k> 510$, we have $-(0.07)k>-(0.076)k+2\geq\theta_3$.  Here note $m_1+m_3\geq k$ by \cite[Theorem 4.4.4]{bcn}. Now $\frac{7}{2}k\geq vk\geq k^2+m_1{\theta_1}^2+m_3{\theta_3}^2$, $\theta_1\geq\frac{1}{4}k$ and $m_1+m_3\geq k$ imply $k\leq510$. This is a contradiction. So, we find that either $m_1<48$ or $k\leq510$. If $m_1<48$, then \cite[Theorem 5.3.2]{bcn} implies $k\leq1127$. In conclusion, we find $k\leq1127$.

   We checked by computer the feasible intersection arrays of the distance-regular graphs with diameter three satisfying  $k_2\leq\frac{3}{2}k$, $\theta_1=\frac{b_1}{2}-1$, $m_1<\frac{1}{2}k$ and $k\leq 1127$, and no  intersection arrays  were feasible.
\end{pf}

\begin{flushleft}
\Large\textbf{Acknowledgments}
\end{flushleft}
We thank E.R. van Dam and S. Bang for their careful reading and comments on an early version of this paper.

\clearpage


\bigskip



\bigskip

\begin{thebibliography}{99}
\bibitem{bhk}
S. Bang, A. Hiraki and J. H. Koolen, Delsarte set graphs with small $c_2$, Graphs Combin. 26 (2010) 147-162

\bibitem{bhk2}
S. Bang, A. Hiraki and J. H. Koolen, Improving diameter bounds for distance-regular graphs, European J. Combin.  27 (2006) 79-89

\bibitem{bbs}
N.L. Biggs, A.G. Boshier and J. Shawe-Taylor, Cubic distance-regular graphs, J. London Math. Soc. 33 (1986) 385-394

\bibitem{bb}
A. Blokhuis, A.E. Brouwer, Determination of the distance-regular graphs without $3$-claws,  Discrete Math. 163 (1997), no. 1-3, 225-227

\bibitem{b}
A.E. Brouwer, Additions and Corrections to Distance-Regular Garphs, http://www.win.tue.nl/$\sim$aeb/drg/index.html

\bibitem{bcn}
A.E. Brouwer, A.M. Cohen and A. Neumaier, Distance-Regular Graphs, Springer-Verlag, Berlin, 1989.

\bibitem{bk}
A.E. Brouwer, J. H. Koolen, The distance-regular graphs of valency four, Journal of Algebraic Combinatorics, 10 (1999) 5-24

\bibitem{bk1}
A.E. Brouwer, J. H. Koolen, The vertex-connectivity of a distance-regular graph, European J. Combin., 30 (2009) 668-673

\bibitem{gr}
C. Godsil and G. Royle, Algebraic Graph Theory, Springer-Verlag, Berlin, 2001.

\bibitem{heam}
W.H. Haemers, Interlacing eigenvalues and graphs, Linear Algebra Appl. 226/228 (1995), 593-616.

\bibitem{kp}
J. H. Koolen and J. Park, Shilla distance-regular graphs, European J. Combin. (2010), doi:10.1016/j.ejc.2010.05.012

\bibitem{kpy}
J. H. Koolen, J. Park and H. Yu, An inequality involving the second largest and smallest eigenvalue of a distance-regular graph, preprint.

\bibitem{mst}
B. Mohar and J. Shawe-Taylor, Distance-biregular graphs with $2$-valent vertices and distance-regular line graphs, J. Combin. Theory Ser. B, 38 (1985), 193-203.

\end{thebibliography}
\end{document}